\documentclass[a4paper]{article}
\usepackage[T1]{fontenc}
\usepackage[utf8]{inputenc}
\usepackage[british]{babel}
\usepackage[top=3cm,bottom=3cm,left=2.5cm,right=2.5cm,heightrounded]{geometry}
\usepackage{amsmath}
\usepackage{amsthm}
\usepackage{amssymb}
\usepackage{braket}
\usepackage{stmaryrd}
\renewcommand{\phi}{\varphi}
\let\originalleft\left
\let\originalright\right
\renewcommand{\left}{\mathopen{}\mathclose\bgroup\originalleft}
\renewcommand{\right}{\aftergroup\egroup\originalright}
\DeclareSymbolFont{bbold}{U}{bbold}{m}{n}
\DeclareSymbolFontAlphabet{\mathbbm}{bbold}
\newcommand{\M}{\mathfrak{M}}
\newcommand{\N}{\mathfrak{N}}
\newcommand{\B}{\mathbb{B}}
\newcommand{\C}{\mathbb{C}}
\newcommand{\0}{\mathbbm{0}}
\newcommand{\1}{\mathbbm{1}}
\newcommand{\Qp}[1]{{\left\llbracket #1 \right\rrbracket}}
\newcommand{\Seq}[1]{\Braket{\, #1 \,}}
\newcommand{\abs}[1]{{\left\lvert #1 \right\rvert}}
\newcommand{\pow}[1]{{\mathcal{P}\left( #1 \right)}}
\newcommand{\fin}[1]{{\left[ #1 \right]^{<\aleph_0}}}
\newcommand{\bp}[2]{{{#1}^{\left[#2\right]}}}
\newcommand{\bu}[3]{{{#1}^{\left[#2\right]}\! / #3}}
\newcommand{\eq}[2]{{\left[ #1 \right]_{#2}}}
\newcommand{\mail}[1]{\href{mailto:#1}{\texttt{#1}}}
\DeclareMathOperator{\dom}{dom}
\DeclareMathOperator{\cf}{cf}
\DeclareMathOperator{\Th}{Th}
\DeclareMathOperator{\sat}{sat}
\theoremstyle{plain}
\newtheorem{theorem}{Theorem}[section]
\newtheorem{proposition}[theorem]{Proposition}
\newtheorem{lemma}[theorem]{Lemma}
\newtheorem{corollary}[theorem]{Corollary}
\newtheorem{conjecture}[theorem]{Conjecture}
\newtheorem{fact}[theorem]{Fact}
\theoremstyle{definition}
\newtheorem{definition}[theorem]{Definition}
\theoremstyle{remark}
\newtheorem{remark}[theorem]{Remark}
\begin{document}
\title{On regular ultrafilters, Boolean ultrapowers, and Keisler's order}
\author{Francesco Parente\\\\ School of Mathematics\\University of East Anglia\\Norwich NR4 7TJ, United Kingdom\\\mail{f.parente@uea.ac.uk}}
\date{}
\maketitle
\begin{abstract}
In this paper we analyse and compare two different notions of regularity for filters on complete Boolean algebras. We also announce two results from a forthcoming paper in preparation, which provide a characterization of Keisler's order in terms of Boolean ultrapowers.
\end{abstract}

\section{Introduction}
Over the last decade, Malliaris and Shelah proved a striking sequence of results in the intersection between model theory and set theory, settled affirmatively the question of whether $\mathfrak{p}=\mathfrak{t}$, and developed surprising connections between classification theory and cardinal characteristics of the continuum. The starting point for their work is the study of Keisler's order, introduced originally in 1967 as a device to compare the complexity of complete theories by looking at regular ultrapowers of their models.

The intuitive idea behind Keisler's order is simple: a theory $T_1$ is “less complicated” than a theory $T_2$ if the ultrapowers of models of $T_1$ are “more likely” to be saturated than the ultrapowers of models of $T_2$. As Malliaris and Shelah~\cite{malliaris:cru} put it, Keisler's order classifies “theories through the lens of ultrafilters”.

It turns out that there is a specific class of ultrafilters which is particularly suitable for this classification work, namely the \emph{regular} ultrafilters.

\begin{definition}[Keisler~\cite{keisler:regular}]\label{definition:regularpow} Let $\kappa$ be an infinite cardinal. A filter $F$ over a set $I$ is \emph{$\kappa$-regular} iff there exists a family $\Set{X_\alpha | \alpha<\kappa}\subseteq F$ such that for every infinite $I\subseteq\kappa$ we have $\bigcap_{\alpha\in I}X_\alpha=\emptyset$.
\end{definition}

The importance of regular ultrafilters lies in the following theorem, which states that whether or not the regular ultrapower of a model of a complete theory is saturated does not depend on the choice of the particular model, but only on the theory itself.

\begin{theorem}[{Keisler~\cite[Corollary 2.1a]{keisler:notsat}}]\label{theorem:kindep} Let $\kappa$ be an infinite cardinal; suppose $U$ is a $\kappa$-regular ultrafilter over a set $I$. If two $L$-structures $\M$ and $\N$ are elementarily equivalent, and $\abs{L}\le\kappa$, then
\[
\M^I\!/U\text{ is }\kappa^+\text{-saturated}\iff\N^I\!/U\text{ is }\kappa^+\text{-saturated}.
\]
\end{theorem}

This theorem then suggests a way of comparing theories according to the saturation of ultrapowers.

\begin{definition}[Keisler~\cite{keisler:notsat}] Let $T_1$ and $T_2$ be complete countable theories and $\kappa$ an infinite cardinal. We define \emph{$T_1\trianglelefteq_\kappa T_2$} iff for every $\kappa$-regular ultrafilter $U$ over $\kappa$ and models $\M_1\models T_1$, $\M_2\models T_2$, if ${\M_2}^\kappa\!/U$ is $\kappa^+$-saturated then ${\M_1}^\kappa\!/U$ is $\kappa^+$-saturated.

Keisler's order is the preorder relation $\trianglelefteq$ defined as follows: $T_1\trianglelefteq T_2$ iff for every infinite $\kappa$, $T_1\trianglelefteq_\kappa T_2$.
\end{definition}

As already mentioned, recent groundbreaking research has shed new light on the structure of Keisler's order. For example, Malliaris and Shelah~\cite{ms:ic} showed that there is an infinite strictly descending chain of theories, and Ulrich~\cite{ulrich:nl} proved the existence of incomparable theories assuming a supercompact cardinal.

Although the definition of Keisler’s order makes use of regular ultrafilters over sets, Malliaris and Shelah~\cite{ms:dl} developed the method of \emph{separation of variables}, which involves a “paradigm shift” towards building ultrafilters on complete Boolean algebras to classify theories. More specifically, suppose we want to construct a regular ultrafilter over $\lambda$ with some specific saturation properties. We build instead an ultrafilter $U$ on a complete Boolean algebra $\B$, together with a suitable surjective homomorphism $j\colon\pow{\lambda}\to\B$, in such a way that $j^{-1}[U]$ will be a regular ultrafilter over $\lambda$ with the desired properties.

On the other hand, in a forthcoming paper by Raghavan and Shelah~\cite{ragsh:bu}, Boolean ultrapowers of forcing iterations are used to force inequalities between cardinal invariants at and above $\omega$. A common thread is the construction of ultrafilters on complete Boolean algebras which will realize or omit some types in the corresponding Boolean ultrapowers.

Motivated by the above results, in this paper we shall focus on the interaction between the combinatorial properties of ultrafilters and the model-theoretic properties of Boolean ultrapowers. More precisely, we ask the following question: what kind of classification can arise when we compare theories according to the saturation of Boolean ultrapowers of their models?

Since we have already seen the crucial role of regular ultrafilters in this context, the first step towards an answer consists in finding the right definition of regularity for ultrafilters on complete Boolean algebras. In fact, two such definitions have already appeared in the literature, both under the name “regular”. This confusion motivates the results of Section~\ref{section:due}, where those two notions are compared and shown not to be equivalent.

In Section~\ref{section:tre}, we analyse the model-theoretic properties of regular ultrafilters in terms of Boolean ultrapowers. In particular, we shall focus on cardinality, cofinality, and universality of Boolean ultrapowers. In each case, one notion of regularity behaves as expected, while the other notion is not well behaved.

In the final section we shall announce two forthcoming results, which provide an answer to our question: Keisler's order can be equivalently characterized in terms of saturation of Boolean ultrapowers. Hence, the model-theoretic properties captured by ultrafilters on complete Boolean algebras are exactly the same as for power-set algebras, thus explaining the success of the “paradigm shift” mentioned above.

\section{Two notions of regularity}\label{section:due}

In this section, we shall present and compare two different definitions of regularity for filters on complete Boolean algebras. As we remarked earlier, both notions have appeared in the literature under the name “regular”. To avoid creating further confusion, we have decided to use the names “regular” and “quasiregular” to distinguish them.

Before we present the first definition, we need to introduce some standard terminology.

\begin{definition} Let $\kappa$ be a cardinal; a Boolean algebra $\B$ is \emph{$\kappa$-c.c.}\ iff every antichain in $\B$ has cardinality less than $\kappa$. The \emph{saturation} of $\B$, denoted by $\sat(\B)$, is the least cardinal $\kappa$ such that $\B$ is $\kappa$-c.c.
\end{definition}

\begin{theorem}[Erd\H{o}s and Tarski~\cite{erdostarski:sat}] If $\B$ is an infinite Boolean algebra, then $\sat(\B)$ is an uncountable regular cardinal.
\end{theorem}

The next remark is straightforward, but will be useful in the proof of Proposition~\ref{proposition:regularkcc}.

\begin{remark}\label{remark:kcc} Suppose $\B$ is a complete Boolean algebra. Then, for every cardinal $\kappa<\sat(\B)$ there exists a maximal antichain $A\subset\B$ with $\abs{A}=\kappa$. To prove this, we note that if $\kappa<\sat(\B)$ then, by definition, $\B$ has an antichain of cardinality $\ge\kappa$. Using Zorn's lemma, we may extend this antichain to a maximal antichain $W$. Since $\kappa\le\abs{W}$, it is possible to partition $W$ into $\kappa$ many non-empty disjoint pieces: $W=\bigcup_{i<\kappa} W_i$. Then clearly $A=\Set{\bigvee W_i | i<\kappa}$ is a maximal antichain in $\B$ such that $\abs{A}=\kappa$.
\end{remark}

We are now ready to state the first main definition, which is due to Shelah~\cite{sh:1064}.

\begin{definition}\label{definition:regular} Let $\B$ be a complete Boolean algebra and $\kappa$ an infinite cardinal. We say that an filter $F$ on $\B$ is \emph{$\kappa$-regular} iff there exist a family $\Set{x_\alpha | \alpha<\kappa}\subseteq F$ and a maximal antichain $A\subset\B$ such that:
\begin{itemize}
\item for every $\alpha<\kappa$ and every $a\in A$, either $a\le x_\alpha$ or $a\wedge x_\alpha=\0$;
\item for every $a\in A$, the set $\Set{\alpha<\kappa | a\le x_\alpha}$ is finite.
\end{itemize}
\end{definition}

It follows immediately from Definition~\ref{definition:regular} that if $F$ is a $\kappa$-regular filter and $\lambda\le\kappa$, then $F$ is also $\lambda$-regular.

\begin{remark}\label{remark:regcc} Let $F$ be a filter on $\B$. If a family $\Set{x_\alpha | \alpha<\kappa}\subseteq F$ and a maximal antichain $A\subset\B$ witness the $\kappa$-regularity of $F$, then necessarily $\kappa\le\abs{A}$. Indeed, for every $\alpha<\kappa$ we can choose some $a_\alpha\in A$ such that $a_\alpha\le x_\alpha$; hence, by $\kappa$-regularity, the map $\alpha\mapsto a_\alpha$ is finite-to-one from $\kappa$ to $A$.
\end{remark}

We now present an existence result for regular ultrafilters; the argument is a simple modification of the construction of Frayne, Morel and Scott~\cite[Theorem 1.17]{fms:reduced}. Also, we remark that more general existence results for regular ultrafilters will appear in Raghavan and Shelah~\cite{ragsh:bu}.

\begin{proposition}\label{proposition:regularkcc} Let $\kappa$ be an infinite cardinal. For a complete Boolean algebra $\B$, the following conditions are equivalent:
\begin{enumerate}
\item\label{proposition:regularkccuno} $\B$ is not $\kappa$-c.c.
\item\label{proposition:regularkccdue} there exists a $\kappa$-regular ultrafilter on $\B$.
\end{enumerate}
\begin{proof} $(\ref{proposition:regularkccdue}\Longrightarrow\ref{proposition:regularkccuno})$ We already know from Remark~\ref{remark:regcc} that if there exists a $\kappa$-regular filter on $\B$, then $\B$ has necessarily an antichain of cardinality $\ge\kappa$.

$(\ref{proposition:regularkccuno}\Longrightarrow\ref{proposition:omegaregulardue})$ Suppose $\B$ is not $\kappa$-c.c. By Remark~\ref{remark:kcc}, we can find a maximal antichain $A=\Set{a_i | i<\kappa}$ in $\B$ such that $\abs{A}=\kappa$. Let us fix an enumeration $\fin{\kappa}=\Set{S_i | i<\kappa}$ and define for every $\alpha<\kappa$
\[
x_\alpha=\bigvee\Set{a_i | \alpha\in S_i}.
\]
Observe that for every $\alpha_1,\dots,\alpha_n<\kappa$ we have
\[
x_{\alpha_1}\wedge\dots\wedge x_{\alpha_n}=\bigvee\Set{a_i | \{\alpha_1,\dots,\alpha_n\}\subseteq S_i}>\0,
\]
hence the family $\Set{x_\alpha | \alpha<\kappa}$ has the finite intersection property, and so it generates an ultrafilter $U$ on $\B$.

To show that $U$ is $\kappa$-regular, we just observe that for each $\alpha<\kappa$ and every $i<\kappa$ we have the two implications
\[
a_i\wedge x_\alpha>\0\implies\alpha\in S_i\implies a_i\le x_\alpha.
\]
From this, it follows immediately that the family $\Set{x_\alpha | \alpha<\kappa}\subseteq U$ and the maximal antichain $A\subset\B$ satisfy the two conditions of Definition~\ref{definition:regular}.
\end{proof}
\end{proposition}

We shall now present a second definition of regularity, which can be found in Koppelberg and Koppelberg~\cite{koppelberg:ultrapower} and Huberich~\cite{hub:reg}. This is arguably the most natural generalization of Definition~\ref{definition:regularpow} to the language of Boolean algebras; however, our choice of terminology “quasiregular” is motivated by the results in Section~\ref{section:tre}, which demonstrate that this natural generalization is in fact not well behaved model theoretically.

\begin{definition}\label{definition:quasiregular} Let $\B$ be a complete Boolean algebra and $\kappa$ an infinite cardinal. We say that a filter $F$ on $\B$ is \emph{$\kappa$-quasiregular} iff there exists a family $\Set{x_\alpha | \alpha<\kappa}\subseteq F$ such that for every infinite $I\subseteq\kappa$ we have $\bigwedge_{\alpha\in I}x_\alpha=\0$.
\end{definition}

Again, it follows from Definition~\ref{definition:quasiregular} that if $F$ is a $\kappa$-quasiregular filter and $\lambda\le\kappa$, then $F$ is also $\lambda$-quasiregular.

The next proposition is straightforward, and justifies our choice of terminology.

\begin{proposition}\label{proposition:regquasireg} Let $\kappa$ be an infinite cardinal. For any complete Boolean algebra $\B$, every $\kappa$-regular filter on $\B$ is also $\kappa$-quasiregular.
\begin{proof} Suppose $F$ is a $\kappa$-regular filter on $\B$; this is witnessed by a family $\Set{x_\alpha | \alpha<\kappa}\subseteq F$ and a maximal antichain $A\subset\B$. We shall prove that for every infinite $I\subseteq\kappa$ we have $\bigwedge_{\alpha\in I}x_\alpha=\0$. To obtain a contradiction, suppose this is not the case. Then $\bigwedge_{\alpha\in I}x_\alpha>\0$, which implies the existence of some $a\in A$ with $a\wedge\bigwedge_{\alpha\in I}x_\alpha>\0$, since $A$ is maximal. Therefore, for every $\alpha\in I$ we have $a\wedge x_\alpha>\0$, which implies $a\le x_\alpha$ by the definition of $\kappa$-regularity. Thus, we have shown that there exists $a\in A$ such that $a\le x_\alpha$ for infinitely many $\alpha$'s, a contradiction.
\end{proof}
\end{proposition}

\begin{lemma} Let $\kappa$ be an infinite cardinal and $\B$ a complete Boolean algebra. If $F$ is a $\kappa$-regular filter on $\B$, then the maximal antichain witnessing its regularity can be chosen to have cardinality $\kappa$.
\begin{proof} Suppose $F$ is a $\kappa$-regular filter on $\B$; this is witnessed by a family $\Set{x_\alpha | \alpha<\kappa}\subseteq F$ and a maximal antichain $A\subset\B$. Consider the following antichain:
\[
W=\Set{\bigwedge_{\alpha\in I}x_\alpha\wedge\bigwedge_{\alpha\notin I}\neg x_\alpha | I\subseteq\kappa}\setminus\{\0\}.
\]

From the definition of $W$, it follows that for every $\alpha<\kappa$ and every $w\in W$, either $w\le x_\alpha$ or $w\wedge x_\alpha=\0$. Furthermore, for every $w\in W$ the set $\Set{\alpha<\kappa | w\le x_\alpha}$ must be finite, otherwise it would contradict the $\kappa$-quasiregularity of $F$ (Proposition~\ref{proposition:regquasireg}).

To see that $W$ is maximal, it suffices to observe for every $a\in A$ there exists a set $I\subseteq\kappa$ such that
\[a\le\bigwedge_{\alpha\in I}x_\alpha\wedge\bigwedge_{\alpha\notin I}\neg x_\alpha.
\]
Hence, $\1=\bigvee A\le\bigvee W$ and so $W$ is maximal.

From Remark~\ref{remark:regcc} we already know that $\kappa\le\abs{W}$. To see that $\abs{W}=\kappa$, observe that whenever $I$ is infinite we must have $\bigwedge_{\alpha\in I}x_\alpha=\0$; therefore we have the equality
\[
W=\Set{\bigwedge_{\alpha\in I}x_\alpha\wedge\bigwedge_{\alpha\notin I}\neg x_\alpha | I\in\fin{\kappa}}\setminus\{\0\},
\]
which gives us $\abs{W}\le\kappa^{<\aleph_0}=\kappa$.
\end{proof}
\end{lemma}

Before the next result, let us recall first that a filter $F$ on a complete Boolean algebra $\B$ is \emph{$\aleph_1$-incomplete} iff there exists a countable subset $X\subseteq F$ such that $\bigwedge X\notin F$. 

\begin{proposition}\label{proposition:omegaregular} Let $\B$ be a complete Boolean algebra. For an ultrafilter $U$ on $\B$, the following conditions are equivalent:
\begin{enumerate}
\item\label{proposition:omegaregularuno} $U$ is $\aleph_0$-regular;
\item\label{proposition:omegaregulardue} $U$ is $\aleph_0$-quasiregular;
\item\label{proposition:omegaregulartre} $U$ is $\aleph_1$-incomplete.
\end{enumerate}
\begin{proof} $(\ref{proposition:omegaregularuno}\Longrightarrow\ref{proposition:omegaregulardue})$ Follows immediately from Proposition~\ref{proposition:regquasireg}.

$(\ref{proposition:omegaregulardue}\Longrightarrow\ref{proposition:omegaregulartre})$ Directly from the definition of $\aleph_0$-quasiregularity, we obtain the existence of some infinite $X\subseteq U$ with $\bigwedge X=\0\notin U$, as desired.

$(\ref{proposition:omegaregulartre}\Longrightarrow\ref{proposition:omegaregularuno})$ Suppose $U$ is $\aleph_1$-incomplete; since $U$ is an ultrafilter, this entails the existence of a countable subset $\Set{x_n | n<\omega}\subseteq U$ such that $\bigwedge_{n<\omega}x_n=\0$. Without loss of generality, we may assume that $x_{n+1}<x_n$ for all $n<\omega$, and $x_0=\1$. Let us define for every $i<\omega$
\[
a_i=x_i\wedge\neg x_{i+1};
\]
it is clear that $A=\Set{a_i | i<\omega}$ is an antichain. Furthermore, $A$ is maximal, because for all $i<\omega$
\[
a_0\vee\dots\vee a_i=x_0\wedge\neg x_{i+1}=\1\wedge\neg x_{i+1}=\neg x_{i+1},
\]
and therefore
\[
\bigvee A=\bigvee_{i<\omega}(a_0\vee\dots\vee a_i)=\bigvee_{i<\omega}\neg x_{i+1}=\neg\bigwedge_{i<\omega}x_{i+1}=\1.
\]

To show that $U$ is $\aleph_0$-regular, it is sufficient to observe that for all $i,n<\omega$ we have the two implications
\[
a_i\wedge x_n>\0\implies n\le i\implies a_i\le x_n.
\]
From this, we deduce that the family $\Set{x_n | n<\omega}$ and the maximal antichain $A$ satisfy the two conditions of Definition~\ref{definition:regular}.
\end{proof}
\end{proposition}

Thus, when $\kappa=\aleph_0$, both regularity properties coincide with $\aleph_1$-incompleteness. When $\kappa$ is arbitrary, an additional distributivity assumption on $\B$ will also make the two properties coincide.

\begin{definition}[Smith and Tarski~\cite{smithtarski:distr}] Let $\kappa$ and $\lambda$ be cardinals. A complete Boolean algebra $\B$ is \emph{$\langle\kappa,\lambda\rangle$-distributive} iff for every function $b\colon\kappa\times\lambda\to\B$ we have
\[
\bigwedge_{\alpha<\kappa}\bigvee_{\beta<\lambda}b(\alpha,\beta)=\bigvee_{f\in{^{\kappa\!}\lambda}}\bigwedge_{\alpha<\kappa}b(\alpha,f(\alpha)).
\]
\end{definition}

\begin{proposition} Let $\kappa$ be an infinite cardinal. If $\B$ is a $\langle\kappa,2\rangle$-distributive complete Boolean algebra, then every $\kappa$-quasiregular filter on $\B$ is $\kappa$-regular.
\begin{proof} Suppose $F$ is a $\kappa$-quasiregular filter on $\B$; by definition, there exists a family $\Set{x_\alpha | \alpha<\kappa}\subseteq F$ such that for every infinite $I\subseteq\kappa$ we have $\bigwedge_{\alpha\in I}x_\alpha=\0$.

Let us define
\[
A=\Set{\bigwedge_{\alpha\in I}x_\alpha\wedge\bigwedge_{\alpha\notin I}\neg x_\alpha | I\subseteq\kappa}\setminus\{\0\};
\]
first of all, it is clear that $A$ is an antichain of $\B$. Furthermore, since each subset $I\subseteq\kappa$ corresponds to its characteristic function $f\colon\kappa\to 2$, we can apply $\langle\kappa,2\rangle$-distributivity to conclude that
\[
\bigvee A=\bigvee\Set{\bigwedge_{\alpha\in I}x_\alpha\wedge\bigwedge_{\alpha\notin I}\neg x_\alpha | I\subseteq\kappa}=\bigwedge_{\alpha<\kappa}(x_\alpha\vee\neg x_\alpha)=\bigwedge_{\alpha<\kappa}\1=\1.
\]
This shows that $A$ is a maximal antichain. By definition of $A$, it follows that for every $\alpha<\kappa$ and every $a\in A$, either $a\le x_\alpha$ or $a\wedge x_\alpha=\0$. Furthermore, for every $a\in A$ the set $\Set{\alpha<\kappa | a\le x_\alpha}$ is finite, otherwise we would contradict the $\kappa$-quasiregularity of $F$. This shows that $F$ is $\kappa$-regular.
\end{proof}
\end{proposition}

\subsection*{The Cohen algebra}

We now focus on a specific complete Boolean algebra which will provide many examples of quasiregular ultrafilters which are not regular.

\begin{definition} For an infinite cardinal $\kappa$, let $\mathbb{P}_\kappa$ be the set of finite partial functions from $\kappa$ to $2$. Given $p,q\in\mathbb{P}_\kappa$, we define $q\le p$ if and only if $p\subseteq q$. Thus, $\langle\mathbb{P}_\kappa,\le\rangle$ is the forcing notion that adjoins $\kappa$ Cohen reals.
\end{definition}

As usual, we say that two conditions $p,q\in\mathbb{P}_\kappa$ are \emph{compatible} iff there exists $r\in\mathbb{P}_\kappa$ such that $r\le p$ and $r\le q$ (otherwise, $p$ and $q$ are \emph{incompatible}).

By a standard result (see Jech~\cite[Corollary 14.12]{JECH}), there exists a unique complete Boolean algebra $\C_\kappa$, usually referred to as the \emph{Cohen algebra}, with a function $e\colon\mathbb{P}_\kappa\to\C_\kappa\setminus\{\0\}$ such that:
\begin{itemize}
\item if $q\le p$ then $e(q)\le e(p)$;
\item $p$ and $q$ are compatible in $\mathbb{P}_\kappa$ if and only if $e(p)\wedge e(q)>\0$;
\item $e[\mathbb{P}_\kappa]$ is dense in $\C_\kappa\setminus\{\0\}$.
\end{itemize}

The following fact is well known and we do not prove it here; a proof can be found, for example, in Jech~\cite{JECH}.

\begin{fact}\label{fact:cohen} The Cohen algebra $\C_\kappa$ is an $\aleph_1$-c.c.\ complete Boolean algebra of cardinality $\kappa^{\aleph_0}$.
\end{fact}

In particular, Proposition~\ref{proposition:regularkcc} implies that no filter on $\C_\kappa$ is $\aleph_1$-regular. On the other hand, the following lemma will provide plenty of $\kappa$-quasiregular ultrafilters on $\C_\kappa$.

\begin{lemma}[Koppelberg and Koppelberg~\cite{koppelberg:ultrapower}]\label{lemma:kopp} On the Cohen algebra $\C_\kappa$, every ultrafilter is $\kappa$-quasiregular.
\begin{proof} For each $\alpha<\kappa$, let us define in $\mathbb{P}_\kappa$
\[
p_{\alpha,0}=\{\langle\alpha,0\rangle\},\quad p_{\alpha,1}=\{\langle\alpha,1\rangle\}.
\]
Firstly, we prove that for every $\alpha<\kappa$
\[
\neg e(p_{\alpha,0})=e(p_{\alpha,1}).
\]
Clearly, since $p_{\alpha,0}$ and $p_{\alpha,1}$ are incompatible, we have $e(p_{\alpha,0})\wedge e(p_{\alpha,1})=\0$. Furthermore, we have $e(p_{\alpha,0})\vee e(p_{\alpha,1})=\1$, because otherwise, by density of $e[\mathbb{P}_\kappa]$, we could find some $p\in\mathbb{P}_\kappa$ such that $e(p)\wedge(e(p_{\alpha,0})\vee e(p_{\alpha,1}))=\0$. But then, as a consequence, $p$ and $p_{\alpha,0}$ are incompatible, which means that $p(\alpha)=1$, but also $p$ and $p_{\alpha,1}$ are incompatible, which means that $p(\alpha)=0$, a contradiction.

Secondly, we prove that whenever $I\subseteq\kappa$ is infinite, we have
\begin{equation}\label{eq:koppuno}
\bigwedge_{\alpha\in I}e(p_{\alpha,0})=\0\quad\text{and}\quad\bigwedge_{\alpha\in I}e(p_{\alpha,1})=\0.
\end{equation}
To obtain a contradiction, suppose that $\bigwedge_{\alpha\in I}e(p_{\alpha,0})>\0$. Then, by density of $e[\mathbb{P}_\kappa]$, there exists some $p\in\mathbb{P}_\kappa$ such that
\begin{equation}\label{eq:kopp}
e(p)\le\bigwedge_{\alpha\in I}e(p_{\alpha,0}).
\end{equation}
We now distinguish two cases, and derive a contradiction in each case. If $p\le p_{\alpha,0}$ for every $\alpha\in I$, then $\bigcup_{\alpha\in I}p_{\alpha,0}\subseteq p$, which is impossible as $p$ is a finite function. On the other hand, if $p\nleq p_{\alpha,0}$ for some $\alpha\in I$, then $\langle\alpha,0\rangle\notin p$, therefore there is some $q\le p$ such that $q(\alpha)=1$. Hence $q\le p_{\alpha,1}$, but then using~\eqref{eq:kopp} we derive
\[
e(q)\le e(p_{\alpha,1})\wedge e(p)=\0,
\]
which is another contradiction. Of course, the same argument also shows that, if $I\subseteq\kappa$ is infinite, then $\bigwedge_{\alpha\in I}e(p_{\alpha,1})=\0$. This completes the proof of~\eqref{eq:koppuno}.

Now let $U$ be any ultrafilter on $\C_\kappa$. If we define
\[
G=\Set{e(p_{\alpha,0}) | \alpha<\kappa},
\]
it follows from what we proved so far that the set
\[
X=(G\cap U)\cup\Set{\neg g | g\in G\setminus U}
\]
is a subset of $U$ with $\abs{X}=\kappa$, such that whenever $Y\subseteq X$ is infinite we have $\bigwedge Y=\0$. Thus, $U$ is $\kappa$-quasiregular.
\end{proof}
\end{lemma}

\subsection*{OK ultrafilters}
While the main focus in this paper is on regular ultrafilters, we conclude this section with a digression on \emph{OK ultrafilters}. Our motivation here is to show that if an ultrafilter is $\aleph_1$-incomplete and $\kappa$-OK, then it is $\kappa$-regular in the sense of Definition~\ref{definition:regular}.

OK ultrafilters were originally defined by Kunen~\cite{kunen:ok} in the context of the topology of $\beta\omega$, the Stone-Čech compactification of the set of natural numbers. Five years later, Dow~\cite{dow:ok} rephrased Kunen's definition in terms of existence of multiplicative functions: this is the definition we present below.

\begin{definition} Let $X$ be any set, $\B$ a Boolean algebra and $f\colon\fin{X}\to\B$.
\begin{itemize}
\item $f$ is \emph{monotonically decreasing} iff for all $S,T\in\fin{X}$, $S\subseteq T$ implies $f(T)\le f(S)$.
\item $f$ is \emph{multiplicative} iff for all $S,T\in\fin{X}$, $f(S\cup T)=f(S)\wedge f(T)$.
\end{itemize}
\end{definition}

\begin{definition}\label{definition:ok} Let $\kappa$ be an infinite cardinal. A filter $F$ on a complete Boolean algebra $\B$ is \emph{$\kappa$-OK} iff for every monotonically decreasing function $f\colon\fin{\kappa}\to F$ such that $\abs{S}=\abs{T}$ implies $f(S)=f(T)$, there exists a multiplicative function $g\colon\fin{\kappa}\to F$ with the property that $g(S)\le f(S)$ for all $S\in\fin{\kappa}$.
\end{definition}

Although it is not completely obvious from Definition~\ref{definition:ok}, it is an easy exercise to verify that if $F$ is a $\kappa$-OK filter and $\lambda\le\kappa$, then $F$ is also $\lambda$-OK.

The model-theoretic relevance of OK ultrafilters lies in a property called \emph{flexibility}, first isolated by Malliaris~\cite{malliaris:flex}. For more details about the connection between OK ultrafilters and Keisler's order we refer the reader to the work of Malliaris and Shelah~\cite{malliaris:cru}.

\begin{theorem}[{Ulrich~\cite[Theorem 5.5]{ulrich:low}}] Let $\kappa$ be an infinite cardinal; suppose $U$ is an $\aleph_1$-incomplete $\kappa$-OK ultrafilter on a complete Boolean algebra $\B$. Then $\B$ is not $\kappa$-c.c.
\end{theorem}

The purpose of the next proposition is to show that Ulrich's argument, which follows the proof of Mansfield~\cite[Theorem 4.1]{MANSFIELD}, can be slightly adapted to obtain a stronger result.

\begin{proposition} Let $\kappa$ be an infinite cardinal; suppose $U$ is an $\aleph_1$-incomplete $\kappa$-OK ultrafilter on a complete Boolean algebra $\B$. Then $U$ is $\kappa$-regular.
\begin{proof} Since $U$ is an $\aleph_1$-incomplete ultrafilter, there exists a countable subset $\Set{a_n | n<\omega}\subseteq U$ such that $\bigwedge_{n<\omega}a_n=\0$. Without loss of generality, we may assume that $a_{n+1}<a_n$ for all $n<\omega$, and $a_0=\1$.

Using this sequence, we can define a monotonically decreasing function as follows:
\[
\begin{split}
f\colon\fin{\kappa} &\longrightarrow U  \\
S &\longmapsto a_\abs{S}
\end{split}
\]
Since $U$ is $\kappa$-OK, we can find a multiplicative function $g\colon\fin{\kappa}\to U$ such that $g(S)\le f(S)$ for all $S\in\fin{\kappa}$. Note that $f(\emptyset)=a_0=\1$, so we may assume without loss of generality that $g(\emptyset)=\1$ as well.

Now, for every $\alpha<\kappa$ define
\[
x_\alpha=g(\{\alpha\}).
\]
For every $n<\omega$, if $\alpha_1,\dots,\alpha_n<\kappa$ are all distinct, then by the multiplicativity of $g$
\[
x_{\alpha_1}\wedge\dots\wedge x_{\alpha_n}=g(\{\alpha_1\})\wedge\dots\wedge g(\{\alpha_n\})=g(\{\alpha_1,\dots,\alpha_n\})\le f(\{\alpha_1,\dots,\alpha_n\})=a_n.
\]
This shows that, whenever $I\subseteq\kappa$ is infinite, we have $\bigwedge_{\alpha\in I}x_\alpha=\0$.

To conclude the proof, we need to find a maximal antichain $A\subset\B$ such that for every $\alpha<\kappa$ and every $a\in A$, either $a\le x_\alpha$ or $a\wedge x_\alpha=\0$. In order to do so, it is sufficient to prove that the set
\[
D=\Set{d\in\B\setminus\{\0\} | \text{for every }\alpha<\kappa\text{, either }d\le x_\alpha\text{ or }d\wedge x_\alpha=\0}
\]
is dense in $\B\setminus\{\0\}$. Then, every maximal antichain $A\subseteq D$ will have the desired property. So let $b\in\B\setminus\{\0\}$; we need to find some $d\in D$ such that $d\le b$. For every $n<\omega$, define
\[
c_n=\bigvee\Set{g(S) | S\in{[\kappa]}^n}.
\]
It is easy to verify that $c_{n+1}\le c_n$ for all $n<\omega$ and that $c_0=g(\emptyset)=\1$. It follows that there exists some $i<\omega$ such that $b\wedge c_i\wedge\neg c_{i+1}>\0$ 
(otherwise, we would have $b\le\bigwedge_{n<\omega}c_n\le\bigwedge_{n<\omega}a_n=\0$, a contradiction). Therefore, by definition of $c_i$, there exists $S\in{[\kappa]}^i$ such that
\[d=b\wedge g(S)\wedge\neg c_{i+1}>\0.
\]
Clearly $d\le b$, so we just need to show that $d\in D$. Let $\alpha<\kappa$; if $\alpha\in S$ then
\[
d=b\wedge g(S)\wedge\neg c_{i+1}\le g(S)\le g(\{\alpha\})=x_\alpha.
\]
Otherwise, if $\alpha\notin S$, then by the multiplicativity of $g$
\[
d\wedge x_\alpha=b\wedge g(S)\wedge g(\{\alpha\})\wedge\neg c_{i+1}=b\wedge g(S\cup\{\alpha\})\wedge\neg c_{i+1}\le b\wedge c_{i+1}\wedge\neg c_{i+1}=\0.
\]
Therefore $d\in D$, as desired.
\end{proof}
\end{proposition}

\section{Model-theoretic properties}\label{section:tre}

In this third section we shall analyse the model-theoretic properties of regular and quasiregular ultrafilters. The natural tool for this analysis is the Boolean ultrapower construction, which dates back to Foster~\cite{foster:bu}. The standard reference for Boolean ultrapowers is Mansfield~\cite{MANSFIELD}; however, since we shall use a slightly different (but equivalent) formulation, the details will be spelled out in the first part of the section. 

Even though we have been working until now with filters on complete Boolean algebras, from now on only ultrafilters will be considered; this is due to the relevance of Theorem~\ref{theorem:los} which we shall be using essentially.

\subsection*{The Boolean ultrapower construction}
Before we present the details of the Boolean ultrapower construction, we need some terminology.

\begin{definition}Let $A$ and $W$ be maximal antichains of a complete Boolean algebra $\B$. We say that $W$ is a \emph{refinement} of $A$ iff for every $w\in W$ there exists $a\in A$ such that $w\le a$. Note that this element $a\in A$ is unique.
\end{definition}

\begin{definition}[Hamkins and Seabold~\cite{HAMSEA}]Let $X$ be any set, $A$ a maximal antichain, and $f\colon A\to X$. If $W$ is a refinement of $A$, the \emph{reduction} of $f$ to $W$ is the function
\[
\begin{split}
(f\mathbin{\downarrow} W)\colon W &\longrightarrow X \\
w &\longmapsto f(a)
\end{split}\ ,
\]
where $a$ is the unique element of $A$ such that $w\le a$.
\end{definition}

\begin{remark} Finitely many maximal antichains $A_1,\dots,A_n$ always admit a common refinement, which is the maximal antichain
\[
\Set{a_1\wedge\dots\wedge a_n | a_i\in A_i}\setminus\{\0\}.
\]
\end{remark}

After these preliminary definitions, we proceed to define a Boolean-valued semantics. The first Definition~\ref{definition:bvsem} deals with the interpretations of the symbols in the language.

\begin{definition}\label{definition:bvsem} Let $\M$ be an $L$-structure and $\B$ a complete Boolean algebra. We define first the set of \emph{names}
\[
\bp{M}{\B}=\Set{\tau\colon A\to M | A\subset\B\text{ is a maximal antichain}}.
\]
\begin{itemize}
\item We now define the Boolean value of the equality symbol: if $\tau,\sigma\in \bp{M}{\B}$, choose a common refinement $W$ of $\dom(\tau)$ and $\dom(\sigma)$, and define
\[
\Qp{\tau=\sigma}^\bp{\M}{\B}=\bigvee\Set{w\in W | (\tau\mathbin{\downarrow} W)(w)=(\sigma\mathbin{\downarrow} W)(w)}.
\]
\item The Boolean values of the symbols in $L$ are defined as follows:
\begin{itemize}
\item if $R\in L$ is an $n$-ary function symbol and $\tau_1,\dots,\tau_n\in \bp{M}{\B}$, choose a common refinement $W$ of $\dom(\tau_1),\dots,\dom(\tau_n)$, and define
\[ \Qp{R(\tau_1,\dots,\tau_n)}^\bp{\M}{\B}=\bigvee\Set{w\in W | \M\models R((\tau_1\mathbin{\downarrow} W)(w),\dots,(\tau_n\mathbin{\downarrow} W)(w))};
\]
\item if $f\in L$ is an $n$-ary function symbol and $\tau_1,\dots,\tau_n,\sigma\in\bp{M}{\B}$, choose a common refinement $W$ of $\dom(\tau_1),\dots,\dom(\tau_n),\dom(\sigma)$, and define
\[ \Qp{f(\tau_1,\dots,\tau_n)=\sigma}^\bp{\M}{\B}=\bigvee\Set{w\in W | \M\models f((\tau_1\mathbin{\downarrow} W)(w),\dots,(\tau_n\mathbin{\downarrow} W)(w))=(\sigma\mathbin{\downarrow} W)(w)};
\]
\item if $c\in L$ is a constant symbol, its interpretation is the name
\[
\begin{split}
c^\bp{\M}{\B}\colon \{\1\} &\longrightarrow M \\
\1 &\longmapsto c^\M
\end{split}\ .
\]
\end{itemize}
\end{itemize}
\end{definition}

Following Mansfield~\cite{MANSFIELD}, the definition of the Boolean values is extended to all formulae in the language, not necessarily atomic: if $\phi(x_1,\dots,x_n)$ is an $L$-formula and $\tau_1,\dots,\tau_n\in\bp{M}{\B}$, the Boolean value
\[
\Qp{\phi(\tau_1,\dots,\tau_n)}^\bp{\M}{\B}
\]
can be defined recursively. From now on, when there is no danger of confusion, the superscript $\bp{\M}{\B}$ will be omitted.

This Boolean-valued semantics is made explicit in the next proposition, which could be taken as a definition:

\begin{proposition}[{Mansfield~\cite[Theorem 1.1]{MANSFIELD}}] \label{proposition:bpowersemantic} Let $\M$ be an $L$-structure and $\B$ a complete Boolean algebra. Let $\phi(x_1,\dots,x_n)$ be an $L$-formula and $\tau_1,\dots,\tau_n\in \bp{M}{\B}$. If $W$ is any common refinement of $\dom(\tau_1),\dots,\dom(\tau_n)$, then
\[
\Qp{\phi(\tau_1,\dots,\tau_n)}=\bigvee\Set{w\in W | \M\models\phi((\tau_1\mathbin{\downarrow} W)(w),\dots,(\tau_n\mathbin{\downarrow} W)(w))}.
\]
\end{proposition}

The following result, sometimes called “mixing property” in the literature, will be useful in the proof of Proposition~\ref{proposition:regcof}.

\begin{proposition}[{Mansfield~\cite[Theorem 1.3]{MANSFIELD}}]\label{proposition:mixing} Let $\M$ be an $L$-structure and $\B$ a complete Boolean algebra. If $A\subset\B$ is an antichain and $\Set{\tau_a | a\in A}\subseteq\bp{M}{\B}$, then there is $\tau\in\bp{M}{\B}$ such that $a\le\Qp{\tau=\tau_a}$ for all $a\in A$.
\end{proposition}

We are now ready to present the main definition.

\begin{definition} Let $\M$ be an $L$-structure, $\B$ a complete Boolean algebra, and $U$ an ultrafilter on $\B$. The \emph{Boolean ultrapower} of $\M$ by $U$, denoted by $\bu{\M}{\B}{U}$, is the $L$-structure defined as follows:
\begin{itemize}
\item Its domain, denoted by $\bu{M}{\B}{U}$, is the quotient of $\bp{M}{\B}$ by the equivalence relation $\equiv_U$ defined as
\[
\tau\equiv_U\sigma\overset{\mathrm{def}}{\iff} \Qp{\tau=\sigma}\in U.
\]
The $\equiv_U$-equivalence class of a name $\tau\in\bp{M}{\B}$ is denoted by $\eq{\tau}{U}$.
\item The interpretations of the symbols in $L$ are defined in the natural way; for example, if $R\in L$ is an $n$-ary relation symbol, then
\[
R^{\bu{\M}{\B}{U}}=\Set{\langle\eq{\tau_1}{U},\dots,\eq{\tau_n}{U}\rangle\in {^{n\!}\bigl(\bu{M}{\B}{U}\bigr)} | \Qp{R(\tau_1,\dots,\tau_n)}\in U},
\]
and similarly for function and constant symbols.
\end{itemize}
\end{definition}

\begin{remark} Suppose $U$ is an ultrafilter over a set $I$. Then, for every structure $\M$
\[\bu{\M}{\pow{I}}{U}\cong\M^I\!/U;
\]
hence, Boolean ultrapowers are indeed a generalization of ultrapowers.
\end{remark}

The following is the analogue for Boolean ultrapowers of a well-known theorem of \L oś~\cite{los:theorem}.

\begin{theorem}[{Mansfield~\cite[Theorem 1.5]{MANSFIELD}}]\label{theorem:los} Let $\M$ be an $L$-structure, $\B$ a complete Boolean algebra, and $U$ an ultrafilter on $\B$. For every $L$-formula $\phi(x_1,\dots,x_n)$ and $\tau_1,\dots,\tau_n\in\bp{M}{\B}$ we have
\[
\bu{\M}{\B}{U}\models\phi(\eq{\tau_1}{U},\dots,\eq{\tau_n}{U})\iff\Qp{\phi(\tau_1,\dots,\tau_n)}\in U.
\]
\end{theorem}

In particular, if for every $m\in M$ we define the name
\[
\begin{split}
\check{m}\colon \{\1\} &\longrightarrow M \\
\1 &\longmapsto m
\end{split}\ ,
\]
then we obtain the following corollary.
\begin{corollary}\label{corollary:los} Let $\M$ be an $L$-structure, $\B$ a complete Boolean algebra, and $U$ an ultrafilter on $\B$. Then the \emph{natural embedding}, defined as 
\[
\begin{split}
d\colon M &\longrightarrow\bu{M}{\B}{U}  \\
m &\longmapsto \eq{\check{m}}{U}
\end{split}\ ,
\]
is an elementary embedding of $\M$ into $\bu{\M}{\B}{U}$.
\end{corollary}

\subsection*{Cardinality}
The problem of determining the possible cardinalities of the ultrapowers of a given structure starts with an simple observation: if $U$ is an ultrafilter over $I$, then for every structure $\M$
\begin{equation}\label{eq:ultcard}
\abs{M}\le\abs{M^I\!/U}\le\abs{M}^\abs{I}.
\end{equation}
Of course, if $U$ is principal then $\abs{M}=\abs{M^I\!/U}$, hence the lower bound in~\eqref{eq:ultcard} can be attained. Therefore, it is natural to ask whether or not the upper bound in~\eqref{eq:ultcard} can be attained for some ultrafilter $U$ over $I$. This question led Frayne, Morel and Scott to consider regular ultrafilters in~\cite{fms:reduced}.

\begin{theorem}[Frayne, Morel and Scott~\cite{fms:reduced}]\label{theorem:fmsreduced} Let $\kappa$ be an infinite cardinal; suppose $U$ is a $\kappa$-regular ultrafilter over a set $I$. For every infinite structure $\M$, we have
\[\abs{M}^\kappa\le\abs{M^I\!/U}.
\]
In particular, if $\abs{I}=\kappa$ then the upper bound $\abs{M}^\abs{I}$ is attained.
\end{theorem}

Motivated by this result, we can ask whether the same is true for Boolean ultrapowers. As we shall see, the parallel of Theorem~\ref{theorem:fmsreduced} is true for regular ultrafilters on complete Boolean algebras, but can fail for quasiregular ultrafilters. First, we need to establish a bound analogous to~\eqref{eq:ultcard}.

\begin{lemma}\label{lemma:bultcard} Let $U$ be an ultrafilter on a complete Boolean algebra $\B$. For every structure $\M$, we have
\begin{equation}\label{eq:bultcard}
\abs{M}\le\abs{\bu{M}{\B}{U}}\le\abs{M}^{<\sat(\B)}+\abs{\B}^{<\sat(\B)}.
\end{equation}
\begin{proof} The inequality $\abs{M}\le\abs{\bu{M}{\B}{U}}$ follows immediately from Corollary~\ref{corollary:los}. On the other hand,
\begin{multline*}
\abs{\bu{M}{\B}{U}}\le \abs{\bp{M}{\B}}=\abs{\Set{\tau\colon A\to M | A\subset\B\text{ is a maximal antichain}}} \\ \le\abs{\bigcup\Set{^X M | X\in\left[\B\right]^{<\sat(\B)}}} = \abs{M}^{<\sat(\B)}+\abs{\B}^{<\sat(\B)},
\end{multline*}
as desired.
\end{proof}
\end{lemma}

We now show that regular ultrafilters produce Boolean ultrapowers of large cardinality; the proof of this result is just a minor modification of the proof of Theorem~\ref{theorem:fmsreduced}.

\begin{proposition} Let $\kappa$ be an infinite cardinal; suppose $U$ is a $\kappa$-regular ultrafilter on a complete Boolean algebra $\B$. For every infinite structure $\M$, we have
\begin{equation}\label{eq:bultregsize}
\abs{M}^\kappa\le\abs{\bu{M}{\B}{U}}.
\end{equation}
In particular, if $\B$ is a $\kappa^+$-c.c.\ Boolean algebra of size $\le 2^\kappa$, then the upper bound in~\eqref{eq:bultcard} is attained.
\begin{proof} Since $\abs{^{<\omega}M}=\abs{M}$, it is sufficient to find an injective function $i\colon{^{\kappa}M}\to\bu{\left(^{<\omega}M\right)}{\B}{U}$. Let the family $\Set{x_\alpha | \alpha<\kappa}\subseteq U$ and the maximal antichain $A\subset\B$ witness the $\kappa$-regularity of $U$. Hence, for every $a\in A$ the set \[S(a)=\Set{\alpha<\kappa | a\le x_\alpha}\] is finite.

Given a function $f\colon\kappa\to M$, we define $\tau_f\colon A\to{^{<\omega}M}$ as follows. Fix $a\in A$; list all the elements of $S(a)$ increasingly as $\alpha_1<\dots<\alpha_n$ and define
\[
\tau_f(a)=\langle f(\alpha_1),\dots,f(\alpha_n)\rangle.
\]

We now prove that the function
\[
\begin{split}
i\colon{^{\kappa}M} &\longrightarrow\bu{\left(^{<\omega}M\right)}{\B}{U}  \\
f &\longmapsto \eq{\tau_f}{U}
\end{split}
\]
is injective. Let $f,g\colon\kappa\to M$; if $f\neq g$ then there exists some $\alpha<\kappa$ such that $f(\alpha)\neq g(\alpha)$. For all $a\in A$, if $a\le x_\alpha$ then $\alpha\in S(a)$ and therefore, by construction, $\tau_f(a)\neq\tau_g(a)$. It follows that
\[
\Qp{\tau_f\neq\tau_g}=\bigvee\Set{a\in A | \tau_f(a)\neq\tau_g(a)}\ge\bigvee\Set{a\in A | a\le x_\alpha}=x_\alpha\in U,
\]
hence $\Qp{\tau_f\neq\tau_g}\in U$, as required. This shows that $i\colon{^{\kappa}M}\to\bu{\left(^{<\omega}M\right)}{\B}{U}$ is injective, establishing~\eqref{eq:bultregsize}.

Now, if we assume further that $\B$ is a $\kappa^+$-c.c.\ Boolean algebra of size $\le 2^\kappa$, then for every infinite structure $\M$
\[
\abs{M}^\kappa\le\abs{\bu{M}{\B}{U}}\le\abs{M}^{<\sat(\B)}+\abs{\B}^{<\sat(\B)}\le\abs{M}^\kappa + \left(2^\kappa\right)^\kappa=\abs{M}^\kappa,
\]
hence we have equality throughout.
\end{proof}
\end{proposition}

Using the Cohen algebra, we can find a counterexample for quasiregular ultrafilters.

\begin{proposition} Let $\kappa$ be an uncountable cardinal. Then there exists a complete Boolean algebra $\B$ and a $\kappa$-quasiregular ultrafilter $U$ on $\B$ such that, for some infinite structure $\M$,
\[
\abs{\bu{M}{\B}{U}}<\abs{M}^\kappa.
\]
\begin{proof} Let $\C_\kappa$ be the Cohen algebra and let $U$ be an ultrafilter on $\C_\kappa$; we know that $U$ is $\kappa$-quasiregular by Lemma~\ref{lemma:kopp}.

Let $\lambda$ be a cardinal such that
\begin{equation}\label{eq:lqr}
\kappa\le\lambda^{\aleph_0}<\lambda^{\kappa}.
\end{equation}
Note that it is always possible to find such a cardinal: for example, if $\lambda\ge\kappa$ is a strong limit cardinal with $\cf(\lambda)=\aleph_1$, then $\lambda$ satisfies~\eqref{eq:lqr}. Now, if $\M$ is a structure with $\abs{M}=\lambda$ then Lemma~\ref{lemma:bultcard} gives us
\[
\abs{\bu{M}{\C_\kappa}{U}}\le\lambda^{<\aleph_1}+\left(\kappa^{\aleph_0}\right)^{<\aleph_1}=\lambda^{\aleph_0}+\kappa^{\aleph_0}=\lambda^{\aleph_0}<\lambda^\kappa,
\]
as desired.
\end{proof}
\end{proposition}

\subsection*{Cofinality}

An important feature of regular ultrafilters is that they produce ultrapowers of large cofinality. We shall now investigate whether the same is true in the context of complete Boolean algebras and Boolean ultrapowers. Again, our results show that regular ultrafilters behave as expected, while quasiregular ultrafilters are not well behaved.

\begin{proposition}\label{proposition:bendaregcof} Let $\kappa$ be an infinite cardinal; suppose $U$ is a $\kappa$-regular ultrafilter over a set $I$. For every infinite cardinal $\lambda$, the ultrapower $\langle\lambda,<\rangle^I\!/U$ has cofinality $>\kappa$.
\end{proposition}

The above result can be found in Benda and Ketonen~\cite[Theorem 1.3]{bendaket:reg}, where it is referred to as a “standard fact”. It appears also in Koppelberg~\cite[Lemma 2]{bkopp:ult}.

By adapting the usual proof of Proposition~\ref{proposition:bendaregcof}, and using the mixing property of Proposition~\ref{proposition:mixing}, we can establish the corresponding result for Boolean ultrapowers.

\begin{proposition}\label{proposition:regcof} Let $\kappa$ be an infinite cardinal; suppose $U$ is a $\kappa$-regular ultrafilter on a complete Boolean algebra $\B$. For every infinite cardinal $\lambda$, the Boolean ultrapower $\bu{\langle\lambda,<\rangle}{\B}{U}$ has cofinality $>\kappa$.
\begin{proof} Let the family $\Set{x_\alpha | \alpha<\kappa}\subseteq U$ and the maximal antichain $A\subset\B$ witness the $\kappa$-regularity of $U$. In particular, this means that for every $a\in A$ the set \[S(a)=\Set{\alpha<\kappa | a\le x_\alpha}\] is finite.

Given any $\Set{\tau_\alpha | \alpha<\kappa}\subset\bp{\lambda}{\B}$, we show that the sequence $\Set{\eq{\tau_\alpha}{U} | \alpha<\kappa}$ is not cofinal in $\bu{\langle\lambda,<\rangle}{\B}{U}$ by finding some $\sigma\in\bp{\lambda}{\B}$ such that $\Qp{\tau_\alpha\le\sigma}\in U$ for all $\alpha<\kappa$.

For every $a\in A$ we wish to define a name $\sigma_a\in\bp{\lambda}{\B}$ such that
\begin{equation}\label{eq:regcof}
\bigwedge_{\alpha\in S(a)}\Qp{\tau_\alpha\le\sigma_a}=\1.
\end{equation}
To do so, consider the finitely many names $\Set{\tau_\alpha | \alpha\in S(a)}$. Bring their domains to a common refinement $W_a$ and define $\sigma_a\colon W_a\to\lambda$ as follows: for all $w\in W_a$
\[
\sigma_a(w)=\max\Set{(\tau_\alpha\mathbin{\downarrow} W_a)(w) | \alpha\in S(a)}.
\]
Clearly $\sigma_a$ will bound each $\tau_\alpha$, for $\alpha\in S(a)$, with Boolean value $\1$, and so~\eqref{eq:regcof} is proved.

Now, use Proposition~\ref{proposition:mixing} to obtain a name $\sigma\in\bp{\lambda}{\B}$ such that $a\le\Qp{\sigma=\sigma_a}$ for each $a\in A$. Since $x_\alpha\in U$ for every $\alpha<\kappa$, to complete the proof it is sufficient to show that
\begin{equation}\label{eq:cofreg}
x_\alpha\le\Qp{\tau_\alpha\le\sigma}.
\end{equation}
For all $a\in A$, if $a\le x_\alpha$ then $\alpha\in S(a)$, hence $\Qp{\tau_\alpha\le\sigma_a}=\1$ and
\[
a\le\Qp{\sigma=\sigma_a}=\Qp{\sigma=\sigma_a}\wedge\1=\Qp{\sigma=\sigma_a}\wedge\Qp{\tau_\alpha\le\sigma_a}\le\Qp{\tau_\alpha\le\sigma}.
\]
Thus we have shown that, for all $a\in A$,  if $a\le x_\alpha$ then $a\le\Qp{\tau_\alpha\le\sigma}$. Now~\eqref{eq:cofreg} follows: for every $\alpha<\kappa$
\[
\Qp{\tau_\alpha\le\sigma}\ge\bigvee\Set{a\in A | a\le\Qp{\tau_\alpha\le\sigma}}\ge\bigvee\Set{a\in A | a\le x_\alpha}=x_\alpha\in U,
\]
thus showing that $\Qp{\tau_\alpha\le\sigma}\in U$.
\end{proof}
\end{proposition}

Since the cofinality of an ordered set is not greater than its cardinality, from the estimate of Lemma~\ref{lemma:bultcard} we already obtain a counterexample for quasiregular ultrafilters. To see this, let $\kappa$ be a cardinal such that $\kappa^{\aleph_0}=\kappa$. If $U$ is any ultrafilter on $\C_\kappa$, then $U$ is $\kappa$-quasiregular, however by Lemma~\ref{lemma:bultcard}
\[
\cf\left(\bu{\langle\kappa,<\rangle}{\C_\kappa}{U}\right)\le\abs{\bu{\kappa}{\C_\kappa}{U}}\le\kappa^{<\aleph_1}+\left(\kappa^{\aleph_0}\right)^{<\aleph_1}=\kappa^{\aleph_0}=\kappa.
\]

Actually, we can prove a more general result.

\begin{proposition} Let $\kappa$ be a regular uncountable cardinal and $\B$ a $\kappa$-c.c.\ complete Boolean algebra. For every ultrafilter $U$ on $\B$, the Boolean ultrapower $\bu{\langle\kappa,<\rangle}{\B}{U}$ has cofinality $\kappa$.
\begin{proof} We observe first that for every $\tau\in\bp{\kappa}{\B}$ there exists some $\alpha<\kappa$ such that $\Qp{\tau\le\check{\alpha}}=\1$. Indeed, given a name $\tau$, the $\kappa$-c.c.\ implies that $\abs{\dom(\tau)}<\kappa$. Since $\kappa$ is a regular cardinal, there exists an $\alpha<\kappa$ such that $\tau(a)\le\alpha$ for all $a\in\dom(\tau)$, as required.

Consequently, the natural embedding
\[
\begin{split}
d\colon\kappa &\longrightarrow\bu{\kappa}{\B}{U}  \\
\alpha &\longmapsto \eq{\check{\alpha}}{U}
\end{split}
\]
is strictly increasing and cofinal in $\bu{\langle\kappa,<\rangle}{\B}{U}$. Hence, the cofinality of $\bu{\langle\kappa,<\rangle}{\B}{U}$ is $\kappa$.
\end{proof}
\end{proposition}

We conclude by mentioning a related result for Boolean ultrapowers of $\langle\omega,<\rangle$. 

\begin{proposition}[{Koppelberg and Koppelberg~\cite[Lemma 3]{koppelberg:ultrapower}}]\label{proposition:koppomega} Let $\kappa$ be a regular cardinal with $\kappa^{\aleph_0}=\kappa$. Then there exists an ultrafilter $U$ on $\C_\kappa$ such that
\[
\cf\left(\bu{\langle\omega,<\rangle}{\C_\kappa}{U}\right)=\abs{\bu{\omega}{\C_\kappa}{U}}=\kappa.
\]
\end{proposition}

Starting from Proposition~\ref{proposition:koppomega}, the topic of the possible cardinality and cofinality of a Boolean ultrapower of $\langle\omega,<\rangle$ was further explored by Koppelberg~\cite{bkopp:ult} and Jin and Shelah~\cite{jinsh:ult}.

\subsection*{Universality}

The third model-theoretic property we consider for Boolean ultrapowers is \emph{universality}. Let us recall first the definition of universal structure.

\begin{definition}[Morley and Vaught~\cite{mv:hum}] Let $\lambda$ be a cardinal. An $L$-structure $\M$ is \emph{$\lambda$-universal} iff for every $L$-structure $\N$, if $\abs{N}<\lambda$ and $\N\equiv\M$ then there is an elementary embedding $j\colon\N\to\M$.
\end{definition}

The following characterization of regularity is implicit in Frayne, Morel and Scott~\cite{fms:reduced} and appears explicitly in Keisler~\cite[Theorem 1.5a]{keisler:notsat}.

\begin{theorem}\label{theorem:reguni} Let $\kappa$ be an infinite cardinal; for an ultrafilter $U$ over a set $I$, the following conditions are equivalent:
\begin{enumerate}
\item $U$ is $\kappa$-regular;
\item for every $L$-structure $\M$, with $\abs{L}\le\kappa$, the ultrapower ${\M^I\!/U}$ is $\kappa^+$-universal.
\end{enumerate}
\end{theorem}

Again, we can adapt the proof of Theorem~\ref{theorem:reguni} to establish a similar characterization of regularity for ultrafilters on complete Boolean algebras.

\begin{proposition}\label{proposition:regularuniversal} Let $\B$ be a complete Boolean algebra and $\kappa$ an infinite cardinal. For an ultrafilter $U$ on $\B$, the following conditions are equivalent:
\begin{enumerate}
\item\label{proposition:regularuniversaluno} $U$ is $\kappa$-regular;
\item\label{proposition:regularuniversaldue} for every $L$-structure $\M$, with $\abs{L}\le\kappa$, the Boolean ultrapower $\bu{\M}{\B}{U}$ is $\kappa^+$-universal.
\end{enumerate}
\begin{proof} $(\ref{proposition:regularuniversaluno}\Longrightarrow\ref{proposition:regularuniversaldue})$ Suppose $U$ is $\kappa$-regular; this is witnessed by a family $\Set{x_\alpha | \alpha<\kappa}\subseteq U$ and a maximal antichain $A\subset\B$. In particular, this means that for every $a\in A$ the set
\[
S(a)=\Set{\alpha<\kappa | a\le x_\alpha}
\]
is finite.

Let $\N$ be an $L$-structure such that $\abs{N}\le\kappa$ and $\N\equiv\bu{\M}{\B}{U}$; we need to find an elementary embedding $j\colon\N\to\bu{\M}{\B}{U}$. Let $L(N)=L\cup\Set{c_n | n\in N}$ be the language obtained from $L$ by adding a new constant symbol $c_n$ for each $n\in N$. We may expand $\N$ to $L(N)$ in a natural way: the interpretation of the symbol $c_n$ is simply $n$; this expansion is denoted by $\N_N$. Let $\Th(\N_N)$ be the set of all $L(N)$-sentences $\phi$ such that $\N_N\models\phi$. Since $\abs{\Th(\N_N)}\le\abs{L}+\abs{N}\le\kappa$, we can enumerate this theory as
\[
\Th(\N_N)=\Set{\phi_\alpha | \alpha<\kappa}.
\]

For each $a\in A$, we proceed to define a sequence $\Seq{\tau_n(a) | n\in N}$ of elements of $M$ in the following way: let $\phi$ be the finite conjunction $\bigwedge_{\alpha\in S(a)}\phi_\alpha$. Let $n_1,\dots,n_k$ be the elements of $N$ appearing as parameters in $\phi$, so that we can write it as $\phi(c_{n_1},\dots,c_{n_k})$.

Since $\phi(c_{n_1},\dots,c_{n_k})\in\Th(\N_N)$, clearly we have
\[
\N\models\exists x_1\dots\exists x_k\phi(x_1,\dots,x_k),
\]
where $x_1,\dots, x_k$ are new variables. But $\N\equiv\bu{\M}{\B}{U}\equiv\M$, and therefore
\[
\M\models\exists x_1\dots\exists x_k\phi(x_1,\dots,x_k).
\]
This allows us to choose $\tau_{n_1}(a),\dots,\tau_{n_k}(a)$ in $M$ such that
\[
\M\models\phi(\tau_{n_1}(a),\dots,\tau_{n_k}(a)).
\]
On the other hand, if $n\in N$ does not appear as a parameter in $\phi$, then we are free to define $\tau_n(a)$ arbitrarily. This completes the definition of the sequence $\Seq{\tau_n(a) | n\in N}$.

Note that for every $n\in N$ we have defined a name $\tau_n\colon A\to M$. We claim that the function
\[
\begin{split}
j\colon N &\longrightarrow\bu{M}{\B}{U}  \\
n &\longmapsto \eq{\tau_n}{U}
\end{split}
\]
is an elementary embedding. Given any formula $\phi_\alpha(c_{n_1},\dots,c_{n_k})\in\Th(\N_N)$, we need to show that $\bu{\M}{\B}{U}\models\phi_\alpha\left(\eq{\tau_{n_1}}{U},\dots,\eq{\tau_{n_k}}{U}\right)$. For all $a\in A$, if $a\le x_\alpha$ then $\alpha\in S(a)$ and therefore, by construction, $\M\models\phi_\alpha(\tau_{n_1}(a),\dots,\tau_{n_k}(a))$. It follows that
\[
\Qp{\phi_\alpha(\tau_{n_1},\dots,\tau_{n_k})}=\bigvee\Set{a\in A | \M\models\phi_\alpha(\tau_{n_1}(a),\dots,\tau_{n_k}(a))}\ge\bigvee\Set{a\in A | a\le x_\alpha}=x_\alpha\in U
\]
which implies by Theorem~\ref{theorem:los} that $j$ is an elementary embedding.

$(\ref{proposition:regularuniversaldue}\Longrightarrow\ref{proposition:regularuniversaluno})$ Let $\M=\left\langle\fin{\kappa},{\subseteq},\Seq{\{\alpha\} | \alpha<\kappa}\right\rangle$ be the structure in the language $L$ with a binary relation symbol for the inclusion and $\kappa$ many constant symbols for the singletons $\{\alpha\}\in\fin{\kappa}$, for $\alpha<\kappa$.

We now define a set of $L$-formulae
\[
\Sigma(x)=\Set{\{\alpha\}\subseteq x | \alpha<\kappa},
\]
and we show that $\Sigma(x)$ is realized in $\bu{\M}{\B}{U}$.

Since every finite subset of $\Sigma(x)$ is realized in $\M$, by compactness there exists a model $\N$ of the theory of $\M$ in which $\Sigma(x)$ is realized. Since $\abs{L}=\kappa$, by Löwenheim-Skolem we may assume that $\abs{N}=\kappa$. We have $\N\equiv\M\equiv\bu{\M}{\B}{U}$, and $\bu{\M}{\B}{U}$ is $\kappa^+$-universal by hypothesis, therefore there exists an elementary embedding $j\colon\N\to\bu{\M}{\B}{U}$. So, if $n\in N$ realizes $\Sigma(x)$ in $\N$, then by elementarity $j(n)$ realizes $\Sigma(x)$ in $\bu{\M}{\B}{U}$. This completes the proof that $\Sigma(x)$ is realized in $\bu{\M}{\B}{U}$.

Now, let $\tau\colon A\to\fin{\kappa}$ be such that $\eq{\tau}{U}$ realizes $\Sigma(x)$ in $\bu{\M}{\B}{U}$. For each $\alpha<\kappa$ define
\[
x_\alpha=\bigvee\Set{a\in A | \M\models\{\alpha\}\subseteq\tau(a)},
\]
and note that $x_\alpha\in U$. To show that $U$ is $\kappa$-regular, we just observe that for each $\alpha<\kappa$ and every $a\in A$ we have the two implications
\[
a\wedge x_\alpha>\0\implies\alpha\in \tau(a)\implies a\le x_\alpha.
\]
From this, we conclude that the family $\Set{x_\alpha | \alpha<\kappa}$ and the maximal antichain $A$ satisfy the two conditions of Definition~\ref{definition:regular}.
\end{proof}
\end{proposition}

Koppelberg and Koppelberg~\cite{koppelberg:ultrapower} showed the existence of a $\kappa$-quasiregular ultrafilter $U$ on $\C_\kappa$ such that, for some $L$-structure $\M$, with $\abs{L}=\aleph_1$, the Boolean ultrapower $\bu{\M}{\C_\kappa}{U}$ is not $\aleph_2$-universal.

With Proposition~\ref{proposition:regularuniversal} available to us, we can give a very simple proof of this fact. Let $\kappa$ be any uncountable cardinal and let $U$ be an ultrafilter on $\C_\kappa$. We already know (Lemma~\ref{lemma:kopp}) that $U$ is $\kappa$-quasiregular, however $U$ cannot be $\aleph_1$-regular, due to the $\aleph_1$-c.c. Therefore, by Proposition~\ref{proposition:regularuniversal} there exists some $L$-structure $\M$, with $\abs{L}\le\aleph_1$, such that the Boolean ultrapower $\bu{\M}{\C_\kappa}{U}$ is not $\aleph_2$-universal.

\section{Keisler's order via Boolean ultrapowers}
In this final section we announce two results from a forthcoming paper in preparation~\cite{parente:kobu}. As we discussed in the introduction, regular ultrafilters play an important role in the classification of theories due to the crucial Theorem~\ref{theorem:kindep}. The first result we announce here is a generalization of Theorem~\ref{theorem:kindep} to regular ultrafilters on complete Boolean algebras.

\begin{theorem}\label{theorem:dist} Let $\kappa$ be an infinite cardinal. Suppose $\B$ is a $\langle\kappa,2\rangle$-distributive complete Boolean algebra and $U$ is a $\kappa$-regular ultrafilter on $\B$. If two $L$-structures $\M$ and $\N$ are elementarily equivalent, and $\abs{L}\le\kappa$, then
\[
\bu{\M}{\B}{U}\text{ is }\kappa^+\text{-saturated}\iff\bu{\N}{\B}{U}\text{ is }\kappa^+\text{-saturated}.
\]
\end{theorem}

While all the model-theoretic properties considered in Section~\ref{section:tre} were generalized smoothly to the context of arbitrary complete Boolean algebras, the analogue of Theorem~\ref{theorem:kindep} was established under an additional distributivity assumption on $\B$. We conjecture that, without this assumption, Theorem~\ref{theorem:dist} can be false in general.

\begin{conjecture}\label{conjecture:kobu} There exists a cardinal $\kappa$, a complete Boolean algebra $\B$, and a $\kappa$-regular ultrafilter $U$ on $\B$ such that, for some $L$-structures $\M\equiv\N$, with $\abs{L}\le\kappa$, the Boolean ultrapower $\bu{\M}{\B}{U}$ is $\kappa^+$-saturated but $\bu{\N}{\B}{U}$ is not $\kappa^+$-saturated.
\end{conjecture}

The above Conjecture~\ref{conjecture:kobu} will be addressed, among other things, in our future work in preparation~\cite{parente:kobu}. For the moment, we just observe that this conjecture has an immediate positive answer if we replace “$\kappa$-regular” with “$\kappa$-quasiregular”.

\begin{proposition} \label{proposition:kcc} Let $\kappa$ be a cardinal with $\kappa^{\aleph_0}=\kappa$. There are two elementarily equivalent $\emptyset$-structures $\M\equiv\N$ such that, for every ultrafilter $U$ on $\C_\kappa$, the Boolean ultrapower $\bu{\M}{\C_\kappa}{U}$ is $\kappa^+$-saturated, but $\bu{\N}{\C_\kappa}{U}$ is not $\kappa^+$-saturated. 
\begin{proof} First, observe that an infinite structure $\M$ in the empty language $L=\emptyset$ is $\kappa^+$-saturated if and only if $\kappa<\abs{M}$.

Now, let $\N=\kappa$. By Lemma~\ref{lemma:bultcard}, for each ultrafilter $U$ on $\C_\kappa$ we have
\[
\abs{\bu{\kappa}{\C_\kappa}{U}}\le\kappa^{<\aleph_1}+\left(\kappa^{\aleph_0}\right)^{<\aleph_1}=\kappa^{\aleph_0}=\kappa,
\]
hence $\bu{\N}{\C_\kappa}{U}$ is not $\kappa^+$-saturated.

On the other hand, let $\M$ be a structure of cardinality at least $\kappa^+$ such that $\M\equiv\N$. Then
\[
\kappa<\abs{M}\le\abs{\bu{M}{\C_\kappa}{U}}
\]
which means that $\bu{\M}{\C_\kappa}{U}$ is $\kappa^+$-saturated.
\end{proof}
\end{proposition}

Thus, not only Theorem~\ref{theorem:dist} fails for quasiregular ultrafilters, but also this failure is due trivially to the cardinality of the Boolean ultrapowers and not to their saturation properties.

We now announce the second result, which answers the question we asked in the introduction, namely: what kind of classification can arise when we compare theories according to the saturation of Boolean ultrapowers of their models? 

When trying to define a Boolean-algebraic analogue of Keisler's order and compare it with the usual one, the first obstacle is that, as far as we know, regular ultrafilters on complete Boolean algebras may not satisfy the generalization of Theorem~\ref{theorem:kindep}. In other words, given a complete theory $T$, whether or not the Boolean ultrapower of a model of $T$ is saturated may depend on the choice of a particular model. However, the next definition is designed to work also in this context.

\begin{definition} Let $\lambda$ be a cardinal and $\B$ a complete Boolean algebra. Suppose $U$ is an ultrafilter on $\B$; we say that $U$ \emph{$\lambda$-saturates} a complete theory $T$ iff for every $\lambda$-saturated model $\M\models T$, the Boolean ultrapower $\bu{\M}{\B}{U}$ is $\lambda$-saturated.
\end{definition}

Using the techniques developed by Malliaris and Shelah~\cite{ms:dl} and Shelah~\cite{sh:1064}, we can establish the following characterization, due to appear in~\cite{parente:kobu}.
\begin{theorem} Let $\kappa$ be an infinite cardinal and $T_1, T_2$ complete countable theories. Then the following are equivalent:
\begin{itemize}
\item $T_1\trianglelefteq_\kappa T_2$;
\item for every $\kappa^+$-c.c.\ complete Boolean algebra $\B$ of size $\le 2^\kappa$, and every $\kappa$-regular ultrafilter $U$ on $\B$, if $U$ $\kappa^+$-saturates $T_2$ then $U$ $\kappa^+$-saturates $T_1$.
\end{itemize}
\end{theorem}

In conclusion, this characterization explains the shift towards constructing regular ultrafilters on complete Boolean algebras: indeed, those ultrafilters are able to detect exactly the same properties of theories as ultrafilters on power-set algebras.

\end{document}